\documentclass[12pt]{article}
\usepackage{amsfonts,amssymb,version}

\def\A{{\mathfrak A}}

\def\B{{\mathfrak B}}

\def\cc{{\mathfrak c}}
\def\mm{\mathfrak m}
\def\LL{{\mathfrak L}}

\def\OO{{\mathcal O}}
\def\pp{\mathfrak p}
\def\pp{\mathfrak p}
\def\qq{\mathfrak q}
\def\rr{\mathfrak r}
\def\R{\mathfrak R}

\def\ker{\mathop{\rm ker}\nolimits}
\def\Max{\mathop{\rm Max}\nolimits}
\def\spec{\mathop{\rm Spec}\nolimits}

\let\ov\overline

\def\stm{\refstepcounter{theorem}\paragraph{\thetheorem}}

\makeatletter \def\l@section{\@dottedtocline{1}{0em}{1.2em}} \makeatother

\begin{document}

\pagestyle{empty}

\centerline{\Large\bf How to construct a closed subscheme, or a} 

\centerline{\Large\bf coherent subsheaf, with prescribed germs}

\bigskip

\centerline{\bf Nitin Nitsure} 

\centerline{\small 
School of Mathematics, Tata Institute of Fundamental Research, Mumbai 400 005}

We show that a closed subscheme $Y$ of a given locally noetherian scheme
$X$ can be constructed by prescribing it germs $Y_x$ 
at all points $x\in X$ in a manner
consistent with specialization of points, provided the resulting 
set $\A$ of all associated points of all the germs $Y_x$ 
is locally finite in $X$.

More generally, we prove a similar result (Theorem \ref{coherent subsheaf})
for constructing a coherent subsheaf $F$ of a coherent sheaf $E$ on $X$ 
by prescribing its stalks $F_x$ at
all points $x\in X$ in a manner consistent with specializations of points, 
again provided the set $\A$ of all associated points of all the 
corresponding local quotients $E_x/F_x$ is locally finite in $X$. 

Given coherent sheaves $E$ and $F$ on $X$, and a prescribed family 
of homomorphisms $(\varphi_x: E_x\to F_x)_{x\in X}$ at the level of stalks
consistent with specialization of points, 
we show (Theorem \ref{homomorphism from germs}) how to construct 
a (unique) homomorphism $\varphi: E\to F$ which has germs 
$(\varphi_x)_{x\in X}$. Equivalently, we show that 
there exists a unique global section $s\in \Gamma(X,F)$ 
of any coherent sheaf $F$ which has a prescribed family of germs 
$(s_x \in F_x)_{x\in X}$ which is consistent with specialization of points.

It is not clear how to formulate an analogous result for constructing
a coherent sheaf $E$ of $\OO_X$-modules in terms of prescribed stalks. 
Even when the set $\B$ of all associated points of all the prescribed  
stalks is locally finite in $X$, such a construction need not succeed
as we show with an example.

The results of this note are elementary and basic, so it is 
possible that they are known to some experts.
However, they are new at least to the author himself, 
who would like these results to become (better) known. 

Besides elementary knowledge of coherent sheaves on schemes 
(for example, as in Chapter 2, Section 5 of [1]), 
we assume knowledge of the basics of associated primes (for example, 
as in Chapter 10, Section 2 of [2]).

\stm\label{coherent subsheaf}{\bf Theorem.} 
{\it Let $X$ be a locally noetherian scheme, 
let $E$ be a coherent sheaf of $\OO_X$-modules on $X$, and let
$(J(x)\subset E_x)_{x\in X}$ be a family of sub-$\OO_{X,x}$-modules 
of the stalks $E_x$ of $E$.
There exists a coherent subsheaf $F\subset E$ such that 
$F_x = J(x)$ for all $x\in X$ if and only if 
both the following conditions (Consistency and Local Finiteness) 
are satisfied.

Consistency:  
If $x_0, x_1 \in X$ such that $x_1$ lies in the closure of $x_0$,
then 
$$J(x_0) = \OO_{X, x_0}\otimes_{\OO_{X, x_1}}J(x_1).$$

Local Finiteness: The set $\A$ which consists 
of all $x\in X$ such that the maximal ideal $\mm_x\subset \OO_{X,x}$ 
is an associated prime to the $\OO_{X,x}$-module $E_x/J(x)$ is 
locally finite.

When it exists, the subsheaf $F$ is unique.
}

\bigskip

Taking $E = \OO_X$ in the above, we get the following result,
which along with the above theorem explains the title of this note.

\stm{\bf Corollary.} 
{\it Let $X$ be a locally noetherian scheme, and let
$(J(x)\subset \OO_{X,x})_{x\in X}$ be a family of ideals.
There exists a closed subscheme $Y\subset X$ such that 
$I_{Y,x} = J(x)$ for all $x\in X$ if and only if 
both the consistency and local finiteness conditions
of the above theorem are satisfied by the 
given family of submodules $J(x)\subset \OO_{X,x}$. 
\hfill $\square$}

\medskip

The uniqueness of $F\subset E$ in Theorem \ref{coherent subsheaf}
is clear. Hence it is enough to prove the existence of $F$ on 
an affine open cover of $X$. Therefore to prove 
Theorem \ref{coherent subsheaf}, it is enough to prove the following
lemma.

\stm\label{submodule lemma}{\bf Lemma.}
{\it Let $A$ be a noetherian ring, and let $E$ be a finite $A$-module. 
For each prime $\pp \subset A$, let there be given a submodule
$J(\pp) \subset E_{\pp}$.
There exists a submodule $F\subset E$ such that 
$F_{\pp} = J(\pp)$ for all primes $\pp\subset A$ if and only if 
both the following conditions (Consistency and Finiteness) 
are satisfied.

Consistency: 
For all primes $\pp \subset \qq \subset A$, 
we have $J(\pp) = J(\qq)_{\pp A_{\qq}} = A_{\pp}\otimes_{A_{\qq}} J(\qq)$. 

Finiteness: 
The set $\A$ of all primes $\pp\subset A$ 
for which $\pp A_{\pp}$ is an associated prime of the $A_{\pp}$-module
$E_{\pp}/J(\pp)$ is a finite set.

Moreover, when $F$ exists it is unique, and we have $Ass_A(E/F) = \A$.
}

\stm\label{intersections}
{\bf Remark.} If $M$ is a finite module 
over a noetherian ring $A$, then the set 
$Ass(M)$ of all associated primes of $M$ is finite. 
For nested primes $\pp \subset \qq$, note that 
$\pp A_{\qq} \in Ass_{A_{\qq}}(M_{\qq})$ if and only if $\pp \in Ass_A(M)$.
If $0\ne w \in M$, then it is an 
elementary fact that 
there exists some $a\in A$ such that the annihilator 
$ann_A(aw)\subset A$ of $aw$ is a prime
ideal $\pp$. Then $\pp\in Ass(M)$, and $w/1 \ne 0 \in M_{\pp}$.
Hence the following equalities hold, 
where $i_{\pp} : M \to M_{\pp} : w\mapsto w/1$ 
denotes the localization homomorphism.
$$\bigcap _{\pp \in \spec A}\, \ker(i_{\pp}) 
  = \bigcap_{\mm \in \Max(A)} \, \ker(i_{\mm}) 
  = \bigcap _{\pp \in Ass(M)}\, \ker(i_{\pp}) 
  = 0 \subset M.$$
Consequently, if $E$ is a finite $A$-module and if
$F\subset E$ any submodule, then the family of submodules 
$(F_{\pp} \subset E_{\pp})_{\pp \in \spec A}$
satisfies the equalities
$$\bigcap _{\pp \in \spec A}\, i_{\pp}^{-1} F_{\pp}  
    = \bigcap_{\mm \in \Max(A)} \, i_{\mm}^{-1} F_{\pp} 
    = \bigcap _{\pp \in Ass(E/F)}\, i_{\pp}^{-1} F_{\pp} 
    = F \subset E.$$


\pagestyle{myheadings}
\markright{Nitin Nitsure: Closed subschemes with prescribed germs.}


\medskip

{\bf Proof of Lemma \ref{submodule lemma}.} \\
{\bf Local case.} If $(A,\mm)$ is noetherian local,
then the candidate $F = J(\mm)$ satisfies the lemma because 
$\pp\subset \mm$ for all $\pp$, and 
the consistency condition give the equality $J(\pp) = J(\mm)_{\pp}$. 

Note that we do not need to assume
the finiteness condition on $\A$ 
in the local case: it gets satisfied automatically, with 
$\A = Ass_A(E/J(\mm))$.
Moreover, by Remark \ref{intersections}, we get the equality
$$\bigcap_{\pp \in \A}\, i_{\pp}^{-1}J(\pp) = J(\mm) \subset E.$$

{\bf General case.} 
When $A$ is not necessarily local, 
our candidate for $F\subset E$ is the submodule defined by 
$$F = \bigcap_{\pp \in \A}\,i_{\pp}^{-1}J(\pp)\subset E$$
which is suggested by Remark \ref{intersections}, and which is 
consistent with what we showed above in the local case.
 
Given any prime $\qq \subset A$, we partition
the set $\A$ as $\A = \LL_{\qq} \coprod \R_{\qq}$,
where
$$\LL_{\qq} = \{ \pp \in \A\,|\, \pp \subset \qq\}
\mbox{ and }
\R_{\qq} = \{ \rr \in \A\,|\, \rr \not\subset \qq\}.$$
Hence we get
$$ F = \left( \bigcap_{\pp \in \LL_{\qq}}\,i_{\pp}^{-1}J(\pp) \right)
\bigcap 
\left( \bigcap_{\rr \in \R_{\qq}}\,i_{\rr}^{-1}J(\rr) \right).$$

The given family $(J(\pp))_{\pp \subset A}$ gives by restriction
a subfamily $(J(\pp))_{\pp \subset \qq}$, which is 
a family of submodules of stalks of the ambient 
finite module $E_{\qq}$
over the local ring $A_{\qq}$. The consistency and finiteness
conditions are satisfied by this family, with the corresponding 
finite set $\A_{\qq}$ of primes in $A_{\qq}$ given by
$\A_{\qq}= Ass_{A_{\qq}}(E_{\qq}/J(\qq)) 
= \{ \pp A_{\qq}\,|\, \pp \in \LL_{\qq}\}$.
By identifying the set $\spec A_{\qq}$ with the subset $\{ \pp \,|\, 
\pp\subset \qq \}$ of $\spec A$, we get a canonical identification
$$\A_{\qq} = \{ \pp A_{\qq}\,|\, \pp \in \LL_{\qq}\} = \LL_{\qq}.$$
Hence the application of the local case of the lemma to the family
$(J(\pp))_{\pp \subset \qq}$ 
over the local ring $A_{\qq}$ gives us the equality
$$J(\qq) = \bigcap_{\pp \in \LL_{\qq}}\,i_{\pp,\qq}^{-1}J(\pp) \subset E_{\qq}$$
where $i_{\pp,\qq}: E_{\qq} \to E_{\pp}$ denotes the localization
homomorphism defined when $\pp\subset \qq$.  
As $i_{\pp,\qq}\circ i_{\qq} = i_{\pp}$ whenever
$\pp\subset \qq \subset A$, applying $i_{\qq}^{-1}$ 
to the above equality gives  
$$F = \left( i_{\qq}^{-1}J(\qq)\right) \bigcap 
\left( \bigcap_{\pp \in \R_{\qq}}\,i_{\pp}^{-1}J(\pp) \right) \subset E.$$

By the above, we have $F\subset i_{\qq}^{-1}J(\qq)$, so 
we get an inclusion $F_{\qq} \subset J(\qq)$. 
To complete the proof of the lemma, we just need to prove
the reverse inclusion $J(\qq)\subset F_{\qq}$. 

For this, let $w\in E$ such that $i_{\qq}(w) \in J(\qq)$.
The consistency condition implies that $i_{\pp}(w) \in J(\pp)$
for all $\pp \in \LL_{\qq}$. 
For each $\rr \in \R_{\qq}$, consider the submodule 
$A_{\rr}\ov{i_{\rr}(w)}\subset  E_{\rr}/J(\rr)$ generated by 
the image $\ov{i_{\rr}(w)}$ of $i_{\rr}(w)$ in $E_{\rr}/J(\rr)$.
We have
$$Ass_{A_{\rr}}(A_{\rr}\ov{i_{\rr}(w)}) \subset Ass_{A_{\rr}}(E_{\rr}/J(\rr))
= \LL_{\rr} \subset \A.$$
As $i_{\pp}(w) \in J(\pp)$ for all $\pp \in \LL_{\qq}$, 
$A_{\pp}\ov{i_{\pp}(w)} = 0$ for  all $\pp \in \LL_{\qq}$.
Hence we get
$$Ass_{A_{\rr}}(A_{\rr}\ov{i_{\rr}(w)}) \subset 
\A - \LL_{\qq} = \R_{\qq}.$$
It follows that for any $\rr \in \R_{\qq}$ we have 
$$\bigcap_{\cc \in \R_{\qq}}\,\cc A_{\rr} \subset 
\sqrt{ann_{A_{\rr}}(\ov{i_{\rr}(w)})}.$$
Now choose an element in $\cc - \qq$ for each $\cc \in \R_{\qq}$, 
and let $a \in A$ be their product, which is defined as 
$\R_{\qq}$ is finite, being a subset of $\A$.  
Then $a \in \cc  - \qq$ for each $\cc \in \R_{\qq}$.
By the above inclusion of ideals, for any sufficiently large integer $N$ 
we have 
$$a^N \in  ann_{A_{\rr}}(\ov{i_{\rr}(w)}) \mbox{ for all } \rr \in \R_{\qq}.$$
Equivalently, we have $a^N i_{\rr}(w) \in J(\rr)$ for all $\rr \in \R_{\qq}$.   
Now recall that $a\not \in \qq$. It therefore 
follows that $i_{\qq}(w) = (a^Nw)/a^N \in E_{\qq}$, 
and the right hand side is in $F_{\qq}$. 

We have thus shown that if $w\in E$ is such that $i_{\qq}(w) \in J(\qq)$,
then $i_{\qq}(w) \in F_{\qq}$. It follows that 
we have the desired inclusion $J(\qq) \subset F_{\qq}$. 
This completes the proof of Lemma \ref{submodule lemma}, 
hence Theorem \ref{coherent subsheaf} is proved. 
\hfill $\square$

\stm\label{first example}
{\bf Example} Let $A = k[t_1,\ldots,t_n]$ be a polynomial
ring in $n$ variables over a field $k$ where $n\ge 1$.
For each prime $\pp \subset A$, let $J(\pp)\subset A_{\pp}$ 
be the ideal defined by 
$$J(\pp) = \left\{ 
\begin{array}{rl}
\pp A_{\pp} & \mbox{ if $\pp$ is maximal in $A$},\\
A_{\pp}     & \mbox{ otherwise}. 
\end{array}\right.$$
This is clearly a consistent family of ideals. 
For each maximal ideal $\mm$, the ideal 
$\mm A_{\mm}$ is associated to $A_{\mm} /\mm A_{\mm}$, so 
$\Max(A) \subset \A$. 
As $\Max(A)$ is infinite, the finiteness 
condition in Lemma \ref{submodule lemma} 
on the set $\A$ of associated primes 
is not satisfied. Hence there does not exist any 
ideal $I\subset A$ which has stalks $J(\pp)\subset A_{\pp}$
for all $\pp$.

\stm\label{homomorphism from germs}{\bf Theorem.}
{\it 
Let $X$ be a locally noetherian scheme, and let $E$ and $F$ be
coherent sheaves of $\OO_X$-modules. Let $(\varphi(x))_{x\in X}$
be a family of $\OO_{X,x}$-linear homomorphisms $\varphi(x): E_x\to F_x$
at the level of stalks, such that whenever $x_0,x_1\in X$ with 
$x_1\in \ov{\{ x_0 \} }$ we have 
$\varphi(x_0) = \varphi(x_1)_{x_0} = 1\otimes \varphi(x_1)$ 
under the natural identification
$$Hom_{\OO_{X,x_0}}(E_{x_0},F_{x_0}) = \OO_{X,x_0}\otimes
_{\OO_{X,x_1}} Hom_{\OO_{X,x_1}}(E_{x_1},F_{x_1}).$$
Then there exists a unique 
$\OO_X$-linear homomorphism $\varphi: E\to F$
such that $\varphi_x = \varphi(x)$ for all $x\in X$.

In particular, taking $E = \OO_X$, given any family $(s(x))_{x\in X}$
of germs $s(x) \in F_x$ such that whenever $x_0,x_1\in X$ with 
$x_1\in \ov{\{ x_0 \} }$ we have $s(x_0) = s(x_1)_{x_0}$,
there exists a unique global section $s\in \Gamma(X,F)$ such that 
$s_x = s(x)$  for all $x\in X$.
}

\medskip
 
{\it Proof.} Uniqueness is clear. As $\underline{Hom}(E,F)$ is a coherent sheaf
and $Hom(E,F) = \Gamma(X,\underline{Hom}(E,F))$, 
to prove the existence part of the theorem, it is enough to consider the 
special case where $E = \OO_X$. Moreover, given uniqueness, it
is enough to prove existence on a Zariski open cover of $X$. Hence 
we can now assume that $X$ is affine. The result now follows for
the next lemma.

\stm\label{germs of section}
{\bf Lemma} {\it Let $A$ be a noetherian ring, 
and $M$ a finite $A$-module. For each prime $\pp\subset A$, 
let there be given an
element $s(\pp)\in M_{\pp}$, such that 
whenever $\pp\subset \qq$ we have  
$s(\pp) = 1\otimes s(\qq) \in 
A_{\pp}\otimes_{A_{\qq}} M_{\qq}  = M_{\pp}$. Then there exists
a unique $v\in M$ such that $i_{\pp}(v) = s(\pp) \in M_{\pp}$ 
for each $\pp$.} 

\medskip

{\it Proof.} For each $\qq \in X = \spec A$, we will produce 
an open neighbourhood $U$ of $\qq$ in $X$ together with a section 
$u \in \Gamma(U, M^{\sim} )$
where $M^{\sim}$ denotes the coherent sheaf on $X = \spec A$ 
defined by the $A$-module $M$, such that for all $\cc \in U$
we have $u_{\cc} = s(\cc) \in M_{\cc}$. As already observed,
such a section $u$ is unique if it exists, and therefore 
such sections on an open cover of $X$ will glue together to define a 
unique global section $v\in M = \Gamma(X, M^{\sim} )$ with the desired
property.  

For any prime $\qq \in A$, we define a partition
$$ Ass_A(M) = \LL_{\qq} \coprod \R_{\qq}$$ 
where
$$\LL_{\qq} = \{ \pp \in \A\,|\, \pp \subset \qq\}
\mbox{ and }\R_{\qq} = \{ \rr \in \A\,|\, \rr \not\subset \qq\}.$$
For the given $\qq$, we choose an element $w\in M$ and 
an element $t \in A - \qq$ such that 
$s(\qq) = w/t \in M_{\qq}$. 
In terms of these elements, let 
$$U = X - Z(t) - \bigcup_{\rr \in \R_{\qq}}\,Z(\rr),$$
where $Z(I)\subset X$ denotes the zero set defined by any $I\subset A$.

As $t$ is invertible in $\Gamma(U, \OO_X)$, 
we get a section
$$u = x/t \in \Gamma(U, M^{\sim} ).$$
We now show that $u$ has the prescribed 
germ $s(\cc)\in M_{\cc}$ at each $\cc \in U$.

If $\pp\subset \cc$ for some $\pp \in Ass_A(M)$, then as
$\cc \in U$ it follows from the definition of $U$ that
$\pp \subset \qq$. In other words, we have 
$$\LL_{\cc} = \{ \pp \in Ass_A(M) \,|\, \pp \subset \cc \}     
\subset \{ \pp \in Ass_A(M) \,|\, \pp \subset \qq \}     
= \LL_{\qq}.$$ 
Now observe that for any prime $\cc$, 
the Remark \ref{intersections} 
applied to the $A_{\cc}$-module $M_{\cc}$ shows that 
the homomorphism
$$\Phi_{\cc}: M_{\cc} \to \bigoplus_{\pp \in \LL_{\cc}} \, M_{\pp}
: v \mapsto (i_{\pp, \cc}(v))_{\pp \in \LL_{\cc}}$$
is injective. By the consistency assumption, 
the image of $s(\cc)$ under $\Phi_{\cc}$ equals 
$(s(\pp))_{\pp \in \LL_{\cc}}$. As $\LL_{\cc} \subset \LL_{\qq}$,
we also have $i_{\pp,\qq}(s(\qq)) = s(\pp)$ for all $\pp \in \LL_{\cc}$.
Hence we have
$$u_{\pp} = s(\pp) \mbox{ for all }\pp \in \LL_{\cc}.$$
This gives the equalities 
$$\Phi_{\cc}( u_{\cc}) 
= (i_{\pp,\cc}(u_{\cc}))_{\pp \in \LL_{\cc}} 
= (u_{\pp})_{\pp \in \LL_{\cc}}
= (s(\pp))_{\pp \in \LL_{\cc}} 
= (i_{\pp,\cc}(s(\cc)))_{\pp \in \LL_{\cc}} 
= \Phi_{\cc}(s(\cc)).$$
By injectivity of $\Phi_{\cc}$ we conclude that 
$u_{\cc} = s(\cc)$, as we wished to prove.
This completes the proof of Lemma \ref{germs of section}, 
hence Theorem \ref{homomorphism from germs} is proved.
\hfill $\square$

\stm{\bf Re-construction of coherent sheaves from germs.} 
Let $X$ be a locally noetherian scheme
and let $E$ be a coherent sheaf of $\OO_X$-modules.
If $x_0, x_1 \in X$ such that $x_1$ lies in the closure of $x_0$,
then we get an induced isomorphism 
$$\rho^E_{x_0, x_1}: \OO_{X, x_0}\otimes_{\OO_{X, x_1}}E_{x_1} \to 
E_{x_0}.$$
These isomorphisms satisfy the co-cycle condition, that is, 
if $x_0,x_1,x_2 \in X$ such that 
$x_1 \in \ov{\{ x_0\} }$ and $x_2 \in \ov{\{ x_1\} }$, then we have
$$\rho^E_{x_0, x_2} = \rho^E_{x_0, x_1}\circ \rho^E_{x_1, x_2}.$$
The set $Ass(E)$ of associated points of $E$ equals the  
set of all $x\in X$ such that the maximal ideal $\mm_x\subset \OO_{X,x}$ 
is an associated prime to the $\OO_{X,x}$-module $E_x$.
Note that $Ass(E)$ is a locally finite set in $X$.
The coherent sheaf $E$ can be reconstructed uniquely from the 
data $(E_x, \rho^E_{x,y})$ associated to it. 
By Theorem \ref{homomorphism from germs}, any homomorphism 
$\varphi: E\to F$ between coherent sheaves on $X$ is the same as 
a consistent family of germs $(\varphi_x)_{x\in X}$. 

The above observation is made precise by the Theorem \ref{fully faithful} 
below. For stating this theorem, we need the following definition.

\stm\label{definition} 
{\bf Definition.} To any locally noetherian scheme $X$, we associate
a category $GermsCoh_X$ defined as follows. 
An {\bf object} of $GermsCoh_X$ is an indexed family
$S = (S(x), \sigma_{x,y})$, where for each $x\in X$, $S(x)$ is a given finite
module over $\OO_{X,x}$, and for each 
$x, y \in X$ such that $y$ lies in the closure of $x$,
we are given an isomorphism 
$\sigma_{x, y}: \OO_{X, x}\otimes_{\OO_{X, y}}S(y) \to 
S(x)$, 
such that the following two conditions (Co-cycle condition and 
Local Finiteness condition) are satisfied.

{\it Co-cycle condition}: If $x_0,x_1,x_2 \in X$ with
$x_1\in \ov{\{ x_0 \} }$ and
$x_2 \in \ov{\{ x_1 \} }$, then we have
$\sigma_{x_0, x_2} = \sigma_{x_0, x_1}\circ \sigma_{x_1, x_2}$.

{\it Local finiteness condition}: For any object 
$S = (S(x), \sigma_{x,y})$, the subset $\B(S)\subset X$ 
which consists of all $x\in X$ such that the maximal 
ideal $\mm_x\subset \OO_{X,x}$ 
is an associated prime to the $\OO_{X,x}$-module $S(x)$, is  
locally finite in $X$. 

If $S = (S(x), \sigma^S_{x,y})$ and $T = (T(x), \sigma^T_{x,y})$ 
are objects of $GermsCoh_X$, then a
{\bf morphism} $S\to T$ in $GermsCoh_X$
is a family $(\varphi(x) : S(x) \to T(x))$
of $\OO_{X,x}$-linear homomorphisms such that 
whenever $y\in \ov{\{ x\} }$ we have
$\varphi(x) \circ \sigma^S_{x,y} = \sigma^T_{x,y}\circ\varphi(y)$.   

We have a natural functor $\pi_X^*: Coh_X \to GermsCoh_X$ 
from the category $Coh_X$ of all coherent sheaves on $X$ to
the category $GermsCoh_X$ defined above, which
sends an object $E$ to the object $(E_x,\rho^E_{x,y})$, 
and sends a morphism $\varphi: E\to F$ to the morphism 
$(\varphi_x: E_x\to F_x): (E_x,\rho^E_{x,y}) \to 
(F_x,\rho^F_{x,y})$. 

In these terms, Theorem \ref{homomorphism from germs} says the following.

\stm\label{fully faithful}
{\bf Theorem.} {\it The functor $\pi_X^*:Coh_X \to GermsCoh_X$ is fully
faithful for any locally noetherian scheme $X$. \hfill $\square$}

\medskip

However, $\pi_X^*:Coh_X \to GermsCoh_X$ is not an 
equivalence of categories in general, as we now show.

\stm{\bf Example where $\pi_X^*$ is not essentially surjective.}
Let $k$ be a field and let $X = {\bf A}^n_k = \spec k[t_1,\ldots, t_n]$ where
$n\ge 1$. Let
$(J(x) \subset \OO_{X,x})$ be the family of ideals defined in 
Example \ref{first example}.
Consider the isomorphisms 
$\sigma_{x,y} : \OO_{X, x}\otimes_{\OO_{X, y}}J(y) \to 
J(x)$ that are induced from the ambient $\OO_X$,
that is, whenever $y\in \ov{\{ x\}}$, the localization homomorphism
$i_{x,y}: \OO_{X,y}\to \OO_{X,x}$ restricts on 
$J(y)\subset \OO_{X,y}$ to induce $\sigma_{x,y}$.
This defines an object $S = (J(x), \sigma_{x,y})$ in $GermsCoh_X$,
as the local finiteness condition is also satisfied:
because each $J(x)$ is non-zero and torsion-free, 
the set $\B(S) \subset X$ of Definition \ref{definition} 
is a singleton set consisting of 
the zero ideal alone.

{\it Note}: In this example, the set $\A$ of all associated points of 
all quotients $\OO_{X,x}/J(x)$ (which was relevant for Theorem
\ref{coherent subsheaf}) is infinite, but all that matters is that
the set $\B(S)$ of all associated points of 
all modules $J(x)$ (which is relevant for 
Definition \ref{definition}) is finite.

Any coherent sheaf $E$ with stalks $J(x)$ will be torsion free, so
it will be locally free on a nonempty open subset of $X$.
However, the set of all closed points is dense in $X$, 
and $J(x) = I_{X,x}$ at any closed point $x \in X = {\bf A}^n_k$, 
which is not a free $\OO_{X,x}$-module whenever $n\ge 2$. 
Hence there is no object $E$ in $Coh_X$  
with $\pi_X^*(E) = (J(x), \sigma_{x,y})$ when $n\ge 2$. 
It is a simple exercise that even for $n=1$, there is no such $E$.

\medskip

We end with the following question suggested by the above example.
Its answer is not known to the author.

\stm{\bf Question.} What is the essential image 
of $\pi_X^* : Coh_X \to GermsCoh_X$?
It would be nice to have some reasonable necessary and sufficient 
condition on a family $(S(x), \sigma_{x,y})$ of germs, 
for it to effectively descend to give a coherent sheaf on $X$. 

\bigskip

\small

{\bf References}

[1] R. Hartshorne: {\it Algebraic Geometry}, Springer 1977.

[2] S. Lang: {\it Algebra}, 3rd edition, Addison-Wesley 1993.

\bigskip

\footnotesize

nitsure@math.tifr.res.in \hfill 15 Dec 2014

\end{document}